\documentclass{article}
\usepackage{amscd,amsfonts,euscript,amsmath,amsxtra,amssymb,amsthm}
\usepackage[english]{babel}
\newtheorem{theorem}{Theorem}%[section]
\newtheorem{lemma}{Lemma}
\newtheorem{proposition}{Proposition}
\newtheorem{definition}{Definition}
\newtheorem{claim}{Claim}
\newtheorem{remark}{Remark}
\newtheorem{convention}{Convention}
\newtheorem{example}{Example}

\theoremstyle{definition}

\theoremstyle{remark}

\numberwithin{equation}{section}

\def\CC{\mathbb C}

\usepackage[all]{xy}
\def\A{{{\mathbb A}}}

\def\E{{{\mathbb E }}}

\def\P{{{\mathbb P }}}
\def\F{{{\mathbb F }}}

\def\L{{{\mathbb L }}}

\def\CC{{{\mathcal C}}}

\def\FF{{{\mathcal F}}}

\def\OO{{{\mathcal O}}}

\def\SS{{{\mathcal S}}}

\def\EExt{{{\mathcal E}xt}}
\def\TTor{{{\mathcal T}or}}
\def\FFitt{{{\mathcal F}itt}}
\def\Hom{{{\rm Hom }}}

\def\Spec{{{\rm Spec \,}}}

\def\Proj{{{\rm Proj \,}}}
\def\Hilb{{{\rm Hilb \,}}}
\def\Univ{{{\rm Univ \,}}}

\def\Quot{{{\rm Quot \,}}}
\def\Supp{{{\rm Supp \,}}}
\def\Pic{{{\rm Pic \,}}}
\def\Stab{{{\rm Stab \,}}}

\def\ker{{{\rm ker \,}}}
\def\rank{{{\rm rank \,}}}

\def\id{{{\rm id }}}

\def\mult{{{\rm mult}}}
\def\length{{{\rm length \,}}}
\def\colength{{{\rm colength \,}}}

\def\tors{{{\rm tors \,}}}
\def\inv{{{\rm inv}}}
\begin{document}
{\sl MSC 14J60, 14D20, 14M27}

 {\sl UDC 512.722+512.723}
\medskip

\begin{center}
{\Large\sc On a new compactification of moduli of vector bundles
on a surface, IV: Nonreduced moduli}\end{center}
\medskip
\begin{center}
Nadezda V. TIMOFEEVA

\smallskip

Yaroslavl' State University

Sovetskaya str. 14, 150000 Yaroslavl', Russia

e-mail: {\it ntimofeeva@list.ru}
\end{center}
\bigskip

\begin{quote}
The construction for nonreduced projective moduli scheme of
semi\-stable admissible pairs is performed. We establish the
relation of this moduli scheme with reduced moduli scheme built up
in the previous article and prove that nonreduced moduli scheme
contains an open subscheme which is isomorphic to moduli scheme of
semistable vector bundles. \\
Bibliography: 10 items.
\end{quote}

\markright{On a new compactification of moduli of vector bundles,
V}

\footnotetext[0]{This work was partially supported by the
Institute of Mathematics "Simion Stoilow"\, of the Romanian
Academy (IMAR) (partnership IMAR -- BITDEFENDER), during author's
visit as invited professor on June -- July  2011.}

\section*{Introduction}
We construct the moduli scheme for  the functor of admissible
semi\-stable pairs $((\widetilde S, \widetilde L), \widetilde E)$
introduced
 in the previous articles \cite{Tim0}-\cite{Tim4} of the author. There we passed to
 the reduction of the intermediate scheme used in the construction.
 This yields in the reduced scheme $\widetilde M$ with no relation
 to (non)reducedness of the functor of moduli. The aim of the
 present article is to give precise sense to this situation and to
 construct (possibly) nonreduced moduli whenever the functor of families is nonreduced.

We start with
\begin{itemize}
\item{smooth irreducible projective algebraic surface $S$ over an
algebraically closed field $k$ of zero characteristic;}
\item{ample invertible $\OO_S$-sheaf $L$ which is fixed and called
{\it polarization} of the surface $S$;}
\item{positive integer $r$ and polynomial $rp_E(n)$ with rational
coefficients to be inter\-pret\-ed as {\it rank} and {\it Hilbert
polynomial} of coherent sheaves on $S$. Hilbert polynomial is
compute with respect to $L$.}
\end{itemize}

We work with {\it admissible semistable pairs} $((\widetilde S,
\widetilde L), \widetilde E).$ Any such pair consists of
\begin{itemize}\item{{\it admissible scheme} of view $\widetilde S=\Proj
\bigoplus (I[t]+(t))^s/(t)^{s+1}$ where $I \subset \OO_S$ is ideal
sheaf of colength equal to $l\le C$. Positive integer $C$ depends
on coefficients of the Hilbert polynomial $rp_E(n)$ (see \S 3).
The case $\widetilde S \cong S$ is also considered. The structure
morphism $\sigma: \widetilde S \to S$ is called {\it canonical
morphism} of the scheme $\widetilde S$. It is identity isomorphism
when $\widetilde S \cong S$;}
\item{ample invertible $\OO_{\widetilde S}$-sheaf $\widetilde L=L^m
\otimes \sigma^{-1}I\cdot \OO_{\widetilde S}.$ It is called the
{\it distinguished polarization} for the scheme $\widetilde S$.
The number $m$ can be chosen common for all $\widetilde S$ (claim
\ref{mmax} \S 3). Fix this $m$ and redenote $L^m$ by $L$ so that
$\widetilde L=L \otimes \sigma^{-1}I \cdot \mathcal{O}_{\widetilde
S}$;}
\item{semistable locally free $\OO_{\widetilde S}$-sheaf
$\widetilde E$ of rank $r$ with Hilbert polynomial $rp_E(n)$ if it
is compute with respect to the distinguished polarization
$\widetilde L$. The sheaf $\widetilde E$ is required to satisfy
the quasi-ideality condition (see (\ref{ei}) below).}
\end{itemize}

It is clear that  admissible scheme of the view $\widetilde
S=$\linebreak $\Proj \bigoplus_{s\ge 0}(I[t]+(t))^s/(t^{s+1})$ can
be naturally represented as a union of irreducible components
$\widetilde S= \bigcup_{i\ge 0}\widetilde S_i$ where the  {\it
main} component $\widetilde S_0=\Proj \bigoplus_{s\ge 0}(I)^s$ is
the blowup of the surface $S$ in the sheaf of ideals $I$ and for
$i>0$ $\widetilde S_i$ are irreducible  {\it additional}
components $\bigcup_{i>0} \widetilde S_i$. As it is shown in
\cite{Tim3}, in this case the additional component can have a
structure of nonreduced scheme. Obviously, admissible scheme
consists of a single component $\widetilde S=\widetilde S_0$ if
and only if it is isomorphic to the initial surface $S$.

The restriction $\sigma_0=\sigma|_{\widetilde S_0}: \widetilde S_0
\to S$ of the canonical morphism  $\sigma$ onto the main component
$\widetilde S_0$ is a blowup morphism.

Following \cite[ch. 2, sect. 2.2]{HL} we recall some definitions.
Let $\CC$ be a category, $\CC^o$ its opposite,
$\CC'={\FF}unct(\CC^o, Sets)$ -- a category of functors to the
category of sets. By Yoneda lemma, the functor $\CC \to \CC':
F\mapsto (\underline F: X\mapsto \Hom_{\CC}(X, F))$ includes $\CC$
as a full subcategory in  $\CC'$.

\begin{definition}\cite[ch. 2, definition 2.2.1]{HL}
The functor ${\mathfrak f} \in {\OO}b\, \CC'$ is {\it
corepres\-ent\-ed by the object} $F \in {\OO}b \,\CC$, if there
exist
 $\CC'$-morphism $\psi : {\mathfrak f} \to
\underline F$ such that any morphism $\psi': {\mathfrak f} \to
\underline F'$ factors through the unique morphism  $\omega:
\underline F \to \underline F'$.
\end{definition}

Let $T$ be a scheme over the field $k$.  Consider families of
semistable pairs
\begin{equation}\label{class} {\mathfrak F}_T= \left\{
\begin{array}{l}\pi: \F \to T,  \;\widetilde \L\in Pic \F ,
\;\forall t\in T \;\widetilde L_t=\widetilde \L|_{\pi^{-1}(t)}\mbox{\rm \; is ample;}\\
(\pi^{-1}(t),\widetilde L_t) \mbox{\rm \;admissible scheme
with distinguished}\\ \mbox{\rm polarisation}; \\
 \widetilde \E - \mbox{\rm locally free } \OO_{\F}-\mbox{\rm
 sheaf};\\
 \chi(\widetilde \E\otimes\widetilde \L^{m})|_{\pi^{-1}(t)})=
 rp_E(m);\\
 ((\pi^{-1}(t), \widetilde L_t), \widetilde \E|_{\pi^{-1}(t)}) - \mbox{\rm (semi)stable pair}
 \end{array} \right\} \end{equation}
 and a functor ${\mathfrak f}: (Schemes_k)^o \to (Sets)$ from the
 category of $k$-schemes to the category of sets. This functor
 assigns to any scheme $T$ the set of equivalence classes $({\mathfrak F}_T/\sim).$

 The equivalence relation $\sim$ is defined as follows.
 Families $((\pi: \F \to T, \widetilde \L),
 \widetilde \E)$ and $((\pi': \F' \to T, \widetilde \L'), \widetilde
 \E')$ of the class $\mathfrak F$ are said to be equivalent (notation:  $((\pi: \F \to T, \widetilde \L),
 \widetilde \E) \sim ((\pi': \F' \to T, \widetilde \L'), \widetilde
 \E')$) if\\
 1) there is an isomorphism $\iota:\F \stackrel{\sim}{\longrightarrow}
 \F'$ such that the diagram \begin{equation*}\xymatrix{\F
 \ar[rd]_{\pi}\ar[rr]_{\sim}^{\iota}&&\F' \ar[ld]^{\pi'}\\
&T }
 \end{equation*} commutes.\\
 2) There are linear bundles $L', L''$ on $T$ such that
 $\iota^{\ast}\widetilde \E' = \widetilde \E \otimes \pi^{\ast}
 L',$ $\iota^{\ast}\widetilde \L' = \widetilde \L \otimes \pi^{\ast}
 L''$.

 \begin{definition} The scheme $\widetilde M$ is a {\it coarse moduli space of
 the functor } $\mathfrak f$ if $\mathfrak f$ is corepresented by
 the scheme $\widetilde M$.
 \end{definition}

 The result of the present article is contained in the following
 theorem.

\begin{theorem} \label{th1} The functor $\mathfrak f$
has a coarse moduli space
 $\widetilde M$ which is a projective Noetherian algebraic scheme of finite
 type. The scheme $\widetilde M $ contains an open sub\-scheme
 $\widetilde M_0$ isomorphic to  open subscheme $\overline M_0$ of
 locally free sheaves
 in the Gieseker --
 Maruyama scheme $\overline M $ corresponding to the same data $r, p_E(n)$.
\end{theorem}

All the reasoning of the present paper is applicable to any
Hilbert polynomial with no relation to the value of discriminant
as well as to the number and geometry of irreducible components in
the corresponding Gieseker -- Maruyama scheme. In general
(reducible) case the theorem provides the existence of the coarse
moduli space for any maximal (under inclusion) irreducible
substack in $\coprod({\mathfrak F}_T/\sim)$ if it contains pairs
$((\pi^{-1}(t), \widetilde L_t), \widetilde \E|_{\pi^{-1}(t)})$
such that $(\pi^{-1}(t), \widetilde L_t)\cong (S,L)$. Such pairs
are called
 {\it $S$-pairs}. We mean under $\widetilde M$
the moduli space of a substack containing semistable $S$-pairs.

The article is organized as follows. \S 1 comprises necessary
additions and corrections for the article \cite{Tim4}. In \S 2 we
recall the well-known definitions and results concerning
(non)reducedness in functorial setting. We prove the proposition
about interplay between the reduction of the functor and of its
moduli scheme. \S 3 is devoted to the construction moduli scheme
$\widetilde M$ in the category of all $k$-schemes. There are two
critical points: bounded\-ness of families to parameterize
(proposition \ref{bnd} \S3) and quasi-projectivity of the
corres\-pond\-ing scheme (proposition \ref{qpr} \S3). \S 4
contains the necessary material about subfunctors and moduli
subspaces and complet\-es the proof of theorem \ref{th1}.

\section{Comments for \cite{Tim4} }

In this section we consider the corrections and additions for the
construction of moduli scheme  $\widetilde M$ done in \cite{Tim4}.
All the results hold; the changes concern some proofs. All
reasonings and results of the article \cite{Tim4} are done for the
functor  $\mathfrak f$ on full subcategory of reduced schemes
$(RSch_k)^o$ in $(Schemes_k)^o$, although this fact was not
reflected explicitly in the text. Respectively, we mention under
$\overline M$ the reduced scheme corresponding to Gieseker --
Maruyama scheme.

Also in the definition for the equivalence of families $((\pi: \F
\to T, \widetilde \L), \widetilde \E)$ и $((\pi': \F' \to T,
\widetilde \L'), \widetilde \E')$ (the requirement 2 following the
diagram (0.1))
 must read as in Introduction of the present article.

In the formulas for the sheaf $L_E$ in \S 2 (preceding the diagram
(2.2)) there are misprints; these formulas must read as
$L_E=\bigwedge^r (E \otimes L^m)= L^m \otimes \det E$ and (below)
$L_E=L^m \otimes c_1$.

Below we give the correction for the requirement of quasi-ideality
for sheaves $\widetilde E$, and corrections in \S 7.

The end of \S 1 must read as follows.

The behavior of vector bundles $\widetilde E$ on additional
components $\widetilde S_i \subset \widetilde S$, $i> 0$ is given
by the following easy computation. Standard exact triple (1.1)
\begin{equation}\label{tri}0\to E \to E^{\vee \vee}\to \varkappa \to 0
\end{equation}
is taken by the functor of direct image  $\sigma_i^{\ast}$ to
the exact sequence
\begin{equation}\label{upsing}\dots \to
\TTor_1^{\sigma_i^{-1}\OO_S}(\sigma^{-1}\varkappa, \OO_{\widetilde
S_i})\to \sigma_i^{\ast}E \to \sigma_i^{\ast}E^{\vee \vee}\to
\sigma_i^{\ast}\varkappa \to 0.\end{equation} In appropriate
neighborhood  $U\subset S$ of the support  $\Supp \varkappa$ the
locally free sheaf $E^{\vee \vee}|_U$ can be replaced by its local
trivialization $\OO_U^{\oplus r}$. Then the exact sequence
(\ref{upsing}) takes the form
\begin{equation*}\label{addex}\dots \to
\TTor_1^{\sigma_i^{-1}\OO_S}(\sigma^{-1}\varkappa, \OO_{\widetilde
S_i})\to \sigma_i^{\ast}E \to \sigma_i^{\ast}\OO_U^{\oplus r}\to
\sigma_i^{\ast}\varkappa \to 0.
\end{equation*}
Consequently, for  $\widetilde
E_i=\sigma^{\ast}E/\tors|_{\widetilde
S_i}=\sigma_i^{\ast}E/\tors_i$ we have
\begin{equation}\label{eei}\dots \to \sigma_i^{\ast}E/\tors_i \to
\sigma_i^{\ast}\OO_U^{\oplus r} \to \sigma_i^{\ast} \varkappa \to
0,\end{equation} where  subsheaf of torsion $\tors_i$ on
(possibly, nonreduced) scheme $\widetilde S_i$ is defined as
before and  $\tors_i=\tors|_{\widetilde S_i}$. Dots on the left
hand side mean the terms which violate exactness. These terms are
not obliged to have positive codimension in $\widetilde S_i$.
\begin{example} Let $\varkappa= k_x$, then $\widetilde S$
consists of two reduced components: $\widetilde S_0$ is a surface
obtained by blowing up of the reduced point $x$ on the surface
$S$, and $\widetilde S_1\cong \P^2$. The morphism $\sigma_1$ is
constant morphism $\sigma_1: \P^2 \to x$. Then $\sigma_1^{\ast}
\varkappa =\sigma_1^{\ast} k_x=\OO_{\P^2}$, and easy counting of
ranks leads to  $\rank \ker (\sigma_1^{\ast}E/\tors_1 \to
\OO_{\widetilde S_1}^{\oplus r})=1$.
\end{example}

Since the sheaf  $\varkappa$ is supported in a finite collection
of points, then the morphism $\OO_U^{\oplus r} \twoheadrightarrow
\varkappa $ can be replaced by the morphism  $\OO_S^{\oplus
r}\twoheadrightarrow \varkappa$.

Let  $q_0: \OO_S^{\oplus r}\twoheadrightarrow \varkappa$ be a
morphism induced by the exact sequence (\ref{tri}). Then we have
\begin{equation}\label{ei} \widetilde E_i=\sigma_i^{\ast}\ker q_0/\tors.
\end{equation}

According to Proposition 1 \cite{Tim4}, for all semistable
coherent sheaves $E$ with fixed Hilbert polynomial $rp_E(m)$ all
sheaves $\widetilde E_i$ on additional components $\widetilde S_i$
can be described by relations of the form (\ref{ei}) for
appropriate $q_0\in \coprod_{l\le c_2} \Quot^l \OO_S^{\oplus r}$.
\vspace{5mm}

The exact sequence  (\ref{eei}) and the relation (\ref{ei})
provide right requirement for quasi-ideality. This requirement
must be of use instead of (1.5) and (1.6) in \cite{Tim4} whenever
it is involved: in the definition of  $S$-(semi)stability
(Definition 6, \S 5), and in proofs of Proposition 10 (\S 6) and
of Lemma 3 (\S 8).

\begin{definition}\label{semistable}  $S$-{\it stable }(respectively, {\it semistable}) {\it pair }
 $((\widetilde S,\widetilde L), \widetilde E)$ is the
following data:
\begin{itemize}
\item{$\widetilde S=\bigcup_{i\ge 0} \widetilde S_i$
admissible scheme, $\sigma: \widetilde S \to S$ canonical
morphism, $\sigma_i: \widetilde S_i \to S$ its restrictions on
components $\widetilde S_i$, $i\ge 0;$}
\item{$\widetilde E$ vector bundle on the scheme $\widetilde S$;}
\item{$\widetilde L \in Pic\, \widetilde S$ distinguished polarization;}
\end{itemize}
such that
\begin{itemize}
\item{$\chi (\widetilde E \otimes \widetilde
L^{n})=rp_E(n);$}
\item{the sheaf $\widetilde E$ on the scheme $\widetilde S$ is
{\it Gieseker-stable} (respectively, {\it Gieseker-semi\-stable}),
i.e. for any proper subsheaf  $\widetilde F \subset \widetilde E$
for  $n\gg 0$
\begin{eqnarray*}
\frac{h^0(\widetilde F\otimes \widetilde L^{n})}{\rank F}&<&
\frac{h^0(\widetilde E\otimes \widetilde L^{n})}{\rank E},
\\ (\mbox{\rm respectively,} \;\;
\frac{h^0(\widetilde F\otimes \widetilde L^{n})}{\rank F}&\leq&
\frac{h^0(\widetilde E\otimes \widetilde L^{n})}{\rank E}\;);
\end{eqnarray*}}
\item{on each additional component  $\widetilde S_i, i>0,$
the sheaf $\widetilde E_i:=\widetilde E|_{\widetilde S_i}$ is {\it
quasi-ideal,} i.e. it has a description of the form ( \ref{ei} )
for some  $q_0\in \bigsqcup_{l\le c_2} \Quot^l \bigoplus^r \OO_S$.
}\end{itemize}
\end{definition}

The end of the proof of Proposition 10 \cite{Tim4} (following the
inclusion of inverse images of trivial sheaves) after the
replacement takes the following form.

There is a commutative triangle
\begin{equation*}\xymatrix{\bigoplus^r \OO_U \ar@{->>}[r]^{q_0}& \varkappa \\
\bigoplus^{r'}\OO_U \ar@{^(->}[u] \ar[ur]^{q'_0} }
\end{equation*}
where the morphism  $q'_0$ is defined as composite map. Applying
the functor of direct image $\sigma^{\ast}$ and restrictions on
each additional component we have
\begin{eqnarray*}\widetilde E_i&=&\sigma_i^{\ast}\ker q_0/\tors_i,\nonumber\\
\widetilde F_i=\sigma^{\ast}F/\tors|_{\widetilde
S_i}&=&\sigma_i^{\ast}\ker q'_0/\tors_i.
\end{eqnarray*}
This completes the proof.

In the proof of Lemma 3 \cite{Tim4} the formula
$$\widetilde E'_i|_{\rm add}=\ker \sigma^{\ast}(\oplus^r
\OO_S\twoheadrightarrow\varkappa)|_{\rm add}
$$
must be replaced by the formula
$$\widetilde E'_i|_{\rm add}=\sigma^{\ast}\ker (\oplus^r
\OO_S\twoheadrightarrow\varkappa)|_{\rm add}/\tors.
$$

The proof of Proposition  13 \cite[\S 7]{Tim4} is wrong: although
the number of isomorphism classes of Artinian quotient algebras
with bounded length for the ring $k[x,y]$ is finite, the set of
equivalence classes of morphisms $\OO_S \twoheadrightarrow
\varkappa$ for $\length \varkappa \le l_0$ is infinite. For
example, let $\varkappa$ be the algebra of dual numbers
corresponding to nonreduced subscheme of length 2 ("infinitesimal
tangent vector"). Choice of various directions of this vector
provides infinite collection of non-equivalent morphisms $\OO_S
\twoheadrightarrow \varkappa$.

Proposition 13 itself, as well as Definition 14 (together with
Remarks 6 and 7) can be (and must be) removed without any loss for
the results of the article.

Proposition 14 must be preceded by the proposition which has in
the text number 15 and must be formulated with respect to
$\diamond$-product $\overline S_{\ast}= \widetilde S \diamond
\widetilde S_{gr}$:
\begin{proposition} If  $\overline \sigma: \overline S_{\ast}\to
S$ is $\diamond$-product in the monoid $\diamondsuit [E]$, then
for any  $E\in [E]$ images $\overline F_i:= \overline
\sigma^{\ast}F_i/\tors$ of sheaves $F_i$ in Jordan -- H\"{o}lder
filtration are locally free and  $\overline \sigma^{\ast} gr_i
(E)= \overline F_i /\overline F_{i-1}.$
\end{proposition}
The proof is preserved literally.

In the proof of Proposition 27 \cite{Tim4} minimal resolution must
be replaced by $\diamond$-product $\overline S_{\ast}=\widetilde S
\diamond \widetilde S_{gr}$.

To construct the scheme  $\widetilde M$ we need a smooth
resolution $\xi: Q' \to Q$ подсхемы $Q \subset \Quot^{rp_E(t)} (V
\otimes L^{-m}),$ where $V= H^0(S, E\otimes L^m)$. Then the
standard resolution of the family of semistable coherent sheaves
with the base  $Q'$ is performed. This leads to a birational
morphism $\phi:\widetilde Q \to Q'$. To proceed further one needs
equivariance of morphisms $\xi$ and $\phi.$ Equivariance in of use
when GIT-quotient $\mu(\widetilde Q)/PGL(V)=\widetilde M$ is
constructed.

Equivariance of smooth resolution can be achieved by replacement
of resolution due to H. Hironaka by the canonical equivariant
resolution done by O. Villamayor and collaborators (cf., for
example, \cite{BEV}). The construction performed in \cite{Tim4},
operates with base schemes with reduced scheme structures. In
particular, the scheme $Q$ is mentioned to be reduced scheme. Then
consider each component of the scheme $\Quot^{rp_E(t)}(V \otimes
L^{-m})$ as reduced scheme, and for  appropriate  $l\gg 0$
consider closed $GL(V)$-equivariant immersion
$$i_l:\Quot^{rp_E(t)}(V \otimes
L^{-m}) \hookrightarrow G(V \otimes H^0(L^{l-m}), rp_E(l)).$$ Then
there is an induced closed equivariant immersion of each component
of the closure  $\overline Q$ of the subscheme  $Q$ in the scheme
$\Quot^{rp_E(t)}(V \otimes L^{-m})$, into the same variety $G(V
\otimes H^0(L^{l-m}), rp_E(l))$. Then we can think of $\overline
Q$ as reduced equidimensional subscheme in $G(V \otimes
H^0(L^{l-m}), rp_E(l))$.

The situation mentioned above allows to apply the result of
\cite[theorem 2.4]{BEV} what provides required $GL(V)$-equivariant
resolution. The algorithm of equivariant resolution consists of a
sequence of blowups in invariant closed subschemes. Let  $X$ be a
scheme acted upon by an algebraic group $G$ $\alpha: G \times X
\to X$, and $\xi: \widehat X \to X$ be an equivariant morphism.
Convince that the action $\alpha$ {\it induces} an action
$\widehat \alpha: G \times \widehat X \to \widehat X$ in the sense
that the commutative square
\begin{equation*}\label{sq}\xymatrix{G \times \widehat X \ar[d]_{(\id_G, \xi)} \ar[r]^{\widehat \alpha}& \widehat X \ar[d]^\xi\\
G \times X \ar[r]^\alpha & X}
\end{equation*} is Cartesian.

Let $T$ be a scheme, $f: T\to G \times X$ and $h: T \to \widehat
X$ its morphisms such that  $\xi \circ h =\alpha \circ f$. We
denote as  $\inv:G\to  G$ the morphism inverting elements of the
group $G$, and as $p_G: G\times X \to G$ the projection on the
first factor. Then the morphism $\varphi: T\to G\times \widehat X$
is uniquely defined by the formula  $\varphi=(p_G\circ f,
\widehat\alpha \circ(\inv \circ p_G\circ f, h))$ what proves the
universality of the square (\ref{sq}).

On the next step of the construction we perform a standard
resolution of the family of semistable coherent sheaves $\E$ on
the $Q$-based trivial family of surfaces $Q\times S$. Note that
choice of $Q$ as smooth resolution of subscheme in the
Grothendieck's scheme $Quot$ is not unique: if there is a
$PGL(V)$-equivariant resolution $Q$ then for further usage any
blowup of $Q$ in $PGL(V)$-equivariant smooth subscheme will fit as
well as $Q$. In this case the construction of standard resolution
described in \cite{Tim1, Tim2} and applied to different schemes
$Q$, leads to different families $\widetilde Q, \pi: \widetilde
\Sigma \to \widetilde Q, \widetilde \E$. Although non-uniquely
definedness of scheme  $\widetilde Q$ does not affect the final
result of the construction.

When the moduli scheme $\widetilde M$ is built up the morphism
 $\mu: \widetilde Q \to \Hilb^{P(t)}G(V,r)$ and
$PGL(V)$-invariant subscheme $\mu(\widetilde Q)$ are involved. The
subscheme yields in  GIT-quotient  $\widetilde M=\mu(\widetilde
Q)/PGL(V)$. By \cite[proposition 18]{Tim4} the subscheme
$\mu(\widetilde Q)$ corresponds to the subset of those closed
points in $\Hilb^{P(t)}G(V,r)$ which are defined by objects of
para\-metriz\-ation  (admissible semistable pairs). This is the
reason why  the subscheme $\mu(\widetilde Q)$ does not depend of
the choice of resolution $Q$ and of the scheme $\widetilde Q$
which is constructed by $Q$.
\begin{definition} The collection $(\widetilde Q,
\pi: \widetilde \Sigma \to \widetilde Q, \widetilde \L, \widetilde
\E)$ where $\pi$ is a flat morphism of schemes, $\widetilde \E$
locally free sheaf flat over $\widetilde Q$, $\widetilde \L$
invertible sheaf which provides distinguished polarizations on
fibres of $\pi$, is called a {\it standard resolution of the
family} $(Q, p: Q\times S \to Q, \L, \E)$ of coherent semistable
sheaves if
\begin{itemize} \item{there are birational morphisms
$\phi: \widetilde Q \to Q$ and $\Phi: \widetilde \Sigma \to
Q\times S$ fitting into the commutative diagram
\begin{equation*} \xymatrix{\widetilde \Sigma \ar[r]^{\Phi}
\ar[d]_{\pi}& Q\times S
\ar[d]^p\\
\widetilde Q \ar[r]^{\phi}&Q}
\end{equation*}}
\item{morphisms $\phi$ and $\Phi$ become isomorphisms when
restricted on subschemes $\widetilde Q_0=\phi^{-1}Q_0$ and
$\widetilde \Sigma_0=\Phi^{-1}(Q_0 \times S)$;}
\item{for each closed point $\widetilde q \in \widetilde Q$
the corresponding member $((\pi^{-1}(\widetilde q), \widetilde
\L|_{\pi^{-1}(\widetilde q)}),$ $\widetilde
\E|_{\pi^{-1}(\widetilde q)})$ of the family  is semistable
admissible pair;}
\item{there is a {\it descent rule}:
$(\Phi_{\ast}\widetilde \E)^{\vee \vee}=\E$.}
\end{itemize}
\end{definition} Let the scheme $Q$ is supplied with an action
$\beta: G\times Q \to Q$ of algebraic group  $G$.
\begin{definition} The standard resolution
$\phi: \widetilde Q \to Q$ is called  {\it equivariant} if there
is an action $ \alpha: G \times \widetilde Q \to \widetilde Q$ of
algebraic group $G$ on the scheme $\widetilde Q$ such that the
diagram
\begin{equation*} \xymatrix{G \times \widetilde Q \ar[d]_{(\id_G,
\phi)} \ar[r]^{ \alpha}& \widetilde Q
\ar[d]^{\phi}\\
G\times Q \ar[r]^{\beta}& Q}\end{equation*} commutes.
\end{definition}

To analyze GIT-stability of points of the subscheme
$\mu(\widetilde Q)$ and to apply  \cite[\S 9]{Tim4} we can replace
(arbitrary) standard resolution $\widetilde Q$ by
$PGL(V)$-equivariant standard resolution.

We point out the equivariant standard resolution  $\widetilde Q'$
such that actions of $G=PGL(V)$ on schemes  $Q$ and
$\mu(\widetilde Q)$ are  $\widetilde Q'$-concordant in the sense
of  \cite[\S 9]{Tim4}. To construct this resolution form a product
of actions $\alpha$ and $\beta$ in the following fashion: $$
\alpha\times \beta: G \times \mu(\widetilde Q) \times Q \to
\mu(\widetilde Q) \times Q : (g, \widetilde q, q) \mapsto
(\alpha(g, \widetilde q), \beta(g, q)).
$$
This is true action of the group $G$ on the product of schemes
$\mu(\widetilde Q) \times Q $. Now consider a locally closed
subscheme ("the diagonal") $(\mu(\phi^{-1}), i):Q_0
\hookrightarrow \mu(\widetilde Q) \times Q $ in this product. The
diagonal is defined on closed points by the correspondence  $q
\mapsto (\mu(\phi^{-1}(q)), q)$. Here $i: Q_0 \hookrightarrow Q$
means the open immersion and we take into account that
$\phi|_{\phi^{-1}(Q_0)}: \widetilde Q_0 \to Q_0$ is an
isomorphism. Note that the image $(\mu(\phi^{-1}), i)(Q_0)$ is
$G$-invariant subscheme with respect to the action $\alpha \times
\beta$. Define a subscheme $\widetilde Q'$ as a closure
$\widetilde Q':=\overline{(\mu(\phi^{-1}), i)(Q_0)}$. It is also
$G$-invariant under the action $\alpha \times \beta$. Then we can
define an action $\widetilde \alpha: G\times \widetilde Q' \to
\widetilde Q'$ by restriction of the product $\alpha \times \beta
$ on the invariant subscheme $\widetilde Q'$: $\widetilde \alpha:=
\alpha \times \beta |_{G\times \widetilde Q'}.$ Note that the
scheme  $\widetilde Q'$ has two surjective morphisms of
projections on factors as defined by the diagram
\begin{equation*}\xymatrix{&\ar[dl]_{p_1}\widetilde Q' \ar@{_(->}[d]\ar[rd]^{p_2}\\
\mu(\widetilde Q)&\ar[l]_{pr_1} \mu(\widetilde Q) \times Q
\ar[r]^{pr_2}& Q}
\end{equation*}

For $\widetilde Q'$-concordance of actions $\alpha$ and $\beta$ it
is sufficient to confirm that both left hand side and right hand
side squares of the commutative diagram
\begin{equation}\label{twosq}\xymatrix{G\times \mu(\widetilde Q)
\ar[d]_{\alpha} &\ar[l]_{(\id_G,p_1)} G\times \widetilde
Q'\ar[d]^{\widetilde \alpha} \ar[r]^{(\id_G,p_2)} &
G\times Q \ar[d]^{\beta}\\
\mu(\widetilde Q)&\ar[l]_{p_1} \widetilde Q' \ar[r]^{p_2}& Q}
\end{equation} are Cartesian.

We perform the proof of the universality for the left hand side
square because the manipulations for the right hand side square
are analogous. Let  $T$ be a scheme, $f: T\to G\times
\mu(\widetilde Q)$ and $h:T\to \widetilde Q'$ are morphisms making
the square
\begin{equation*}\xymatrix{G\times \mu(\widetilde Q) \ar[d]_\alpha & \ar[l]_f T \ar[d]^h\\
\mu(\widetilde Q)&\ar[l]_{p_1} \widetilde Q'}
\end{equation*} to commute. Let  $t\in T$ be a closed point,
$f(t)=(g, \widetilde q) \in G\times \mu(\widetilde Q)$ and
$h(t)=\widetilde q'$ its images under morphisms $f$ and $h$
respectively. Let also  $\alpha(g, \widetilde q)=g\widetilde q$.
By commutativity of the left hand side square in (\ref{twosq})
$p_1(g\widetilde q)=\widetilde q'$. Then $(\id_G, p_1)(g,
\widetilde q)=(g, p_1(\widetilde q))=(g, g^{-1}\widetilde q).$
Define a map $\varphi: T \to G \times \widetilde Q'$ by the
correspondence $t \mapsto (g, g^{-1}\widetilde q').$ It is clear
that the map $\varphi$ is uniquely defined. Denoting by the symbol
$\inv: G \to G$ the inverting morphism in the group $G$ and by the
symbol $p_G: G \times \mu(\widetilde Q) \to G$ projection to the
factor we define a morphism  $\varphi$ by the formula
$$\varphi=(p_G \circ f, \widetilde \alpha \circ ( \inv \circ p_G
\circ f, h)).$$

Author presents her deepest apologies for mistakes and gaps in the
text of the article \cite{Tim4}.

\section{The reduction of moduli functor and its moduli scheme}
First recall some  definitions \cite[ch. 1, sect. 4.5]{EGA}.

Let $X=(X, \OO_X)$ be a scheme and let $\mathcal{N}il_X \subset
\OO_X$ be the nilradical of the structure sheaf $\OO_X$. Then the
quotient sheaf $\OO_X/\mathcal{N}il_X$ will be denoted as
$\OO_{X_{red}}$. The scheme $X_{red}:=(X, \OO_{X_{red}})$ will be
referred to as a \textit{reduction of the scheme} $(X, \OO_X)$. It
is clear that the reduction of the scheme $(X, \OO_X)$ is
homeomorphic to the scheme $X$ as Zariski topological space and
that it is canonically embedded subscheme in $X$ \cite[ch. 1,
corollaire 4.5.2]{EGA}. If $\rho: X_{red}\hookrightarrow X$ be the
morphism of immersion then the corresponding sheaf morphism is the
morphism onto the quotient sheaf $\rho^{\sharp}: \OO_X \to
\OO_{X_{red}}$.

\begin{proposition} \cite[proposition 4.5.10]{EGA} For any scheme
morphism $f: X\to Y$ there is a canonically defined morphism of
the corresponding reduced schemes $f_{red}: X_{red} \to Y_{red}$
making the diagram
\begin{equation*} \xymatrix{X\ar[r]^f&Y\\
X_{red}\ar@{^{(}->}[u] \ar[r]^{f_{red}} &Y_{red}\ar@{^{(}->}[u]}
\end{equation*}
to commute. \end{proposition} This morphism will be called
\textit{the reduction of the morphism} $f$.

Now consider any functor of the following view
\begin{equation*}
 \mathfrak{f}: (Schemes_k)^o\to Sets
\end{equation*}
from the category of schemes over the field $k$ to the category of
sets. The functor of our interest attaches to any scheme $T\in
\mathrm{Ob} Schemes_k$ the set $\mathfrak{F}_T$ of $T$-based
families of objects of some prescribed type. Let $(RSch_k)^o$ be a
subcategory formed by reduced schemes and their morphisms in
$(Schemes_k)^o$. Then we have a natural inclusion (as a full
subcategory) $\mathfrak{i}:(RSch_k)^o \hookrightarrow
(Schemes_k)^o$ and a functor of reduction $\mathfrak{red}:
(Schemes_k)^o \to (RSch_k)^o$ such that the composite
$\mathfrak{red} \circ \mathfrak{i}$ is identity functor on
$(RSch_k)^o$.

We call the restriction
$\mathfrak{f}_{red}:=\mathfrak{f}|_{(RSch_k)^o}$ of the functor
$\mathfrak{f}$ onto subcategory $(RSch_k)^o$ as \textit{reduction
of the functor} $\mathfrak{f}$. It is natural to say that the
functor $\mathfrak f$ is {\it reduced} if it factors as
${\mathfrak f}={\mathfrak f}_{red} \circ {\mathfrak{red}}$.

The same can be done if the corresponding functor of moduli
\begin{equation*}
\mathfrak{f/\sim}: (Schemes_k)^o \to Sets
\end{equation*}
is considered. This means that it attaches to any scheme $T$ the
set of classes of $T$-based families $\mathfrak{F}_T/\sim$ with
respect to some appropriate equivalence relation $\sim$.

\begin{proposition}
Let $\mathfrak{f}$ be corepresented by a scheme $M$. Then its
reduction $\mathfrak{f}_{red}$ is corepresented by the reduction
$M_{red}$ of the scheme $M$.
\end{proposition}
\begin{proof} Since $\mathfrak{f}$ is corepresented by the scheme $ M$
then there exist a morphism $\psi: \mathfrak{f} \to (X\mapsto
\Hom(X,M))$ such that any morphism $\psi': \mathfrak{f} \to (X
\mapsto \mathrm{Hom}(X, F'))$ factors through the unique morphism
$$\omega: (X \mapsto \Hom(X, M)) \to (X\mapsto \Hom(X,F')).$$
Passing to reductions of $\mathfrak{f}, X, F', M$ one gets the
required. \end{proof}

\section{A nonreduced moduli scheme for $\mathfrak f$}
In this section we construct (possibly) nonreduced moduli scheme
for admissible S-semistable pairs. Afterwards we prove that this
scheme is quasi-projective. In the situation of interest this is
enough to conclude that the nonreduced moduli scheme is projective
scheme.

To construct moduli space one encounters two problems: to examine
the bounded\-ness of families of the interest and to apply
geometric invariant theory in his situation. In our case due to
the special choice of distinguished polarizations $\widetilde L$
Hilbert polynomials both for admissible schemes $\widetilde S$ and
for semistable vector bundles $\widetilde E$ remain constant over
the base.

\begin{claim} \label{mmax} There exist a (common for all semistable admissible
pairs) integer $m \gg 0$ such that the invertible sheaf
$\widetilde L= L^m \otimes \sigma^{-1}I \cdot
\mathcal{O}_{\widetilde S}$ is very ample. \end{claim}

Now we explain why some finite (may be big enough) integer $m$ can
be sufficient for all schemes $\widetilde S$. Isomorphy classes of
schemes $\widetilde S= \mathrm{Proj} \bigoplus_{s\ge 0}
(I[t]+(t))^s/(t)^{s+1}$ are enumerated by 0-dimensional subschemes
in $S$ corresponding to the ideals $I$.

It is clear from the construction of the schemes $\widetilde S$
\cite{Tim3} that lengths of these 0-dimension\-al subschemes are
bounded by the function of coefficients of the Hilbert polynomial
of sheaves of our interest: $l\le C$. Indeed, the sheaf of ideals
$I$ is defined as 0-th Fitting ideal sheaf $\FFitt^0
\EExt^2(\varkappa , \OO_S)$ for $\varkappa$ being the quotient
sheaf $E^{\vee \vee}/E$. It is compute by coherent semistable
torsion-free $\OO_S$-sheaf $E$ with Hilbert polynomial $rp_E(n)$
and defines its resolution of singularity. Note that  $\length
\EExt^2(\varkappa, \OO_S)=\length \varkappa$ and the last length
is bounded from above by the second Chern class $\length \varkappa
\le c_2$ of sheaves of interest. Note that there is only finite
collection of isomorphy classes of
 Artinian quotient $k[x,y]_{(x,y)}$-modules of
$k[x,y]_{(x,y)}^{\oplus r}$ of the fixed length. Here
$k[x,y]_{(x,y)}$ means the localization of the ring $k[x,y]$ in
the maximal ideal $(x,y)$. So we come to the conclusion that the
set of colengths $\colength \FFitt^0 \EExt^2(\varkappa, \OO_S)$ is
finite and hence bounded from above.

Then all isomorphy classes of $\widetilde S \not \cong S$ can be
depicted as points of the following scheme of finite type:
$\coprod_{l=1}^{C}\mathrm{Hilb}^l S$

\begin{convention} Fix this $m$ and redenote $L^m$ by $L$ so that $\widetilde
L=L \otimes \sigma^{-1}I \cdot \mathcal{O}_{\widetilde S}$ as was
mentioned in the introduction. Also for any flat family of schemes
of the class $\widetilde S$ with distinguished polarizations $
\widetilde L $, let $\widetilde{\mathbb{L}}$ be an invertible
sheaf which is very ample relatively to $T$. It is assumed tacitly
to give distinguished polarization on each fibre when restricted
onto this fibre.\end{convention}

Let $X$ be a projective scheme over a field $k$ and let $\OO(1)$
be a very ample line bundle.
\begin{definition}\cite[Definition 1.7.5]{HL} A family of
isomorphism classes of coherent sheaves on $X$ is {\it bounded} if
there is a $k$-scheme $S$ of finite type and a coherent
$\OO_{S\times X}$-sheaf $F$ such that the given family is
contained in the set $$\{F|_{\Spec k_s \times X}|s \mbox{ a closed
point in }S\}.$$
\end{definition}

\begin{theorem}\label{bndth} \cite[Theorem 3.3.7]{HL} Let $f:  X \to S$ be a projective morphism of
schemes of finite type over $k$ and let $\OO_X(1)$ be an $f$-ample
line bundle. Let $P$ be a polynomial of degree $d$, and let
$\mu_0$ be a rational number. Then the family of purely
$d$-dimensional sheaves on the fibres of $f$ with Hilbert
polynomial $P$ and maximal slope $\widehat \mu_{max} \le \mu_0$ is
bounded. In particular, the family of semistable sheaves with
Hilbert polynomial $P$ is bounded.
\end{theorem}

\begin{proposition}\label{bnd} Family of admissible semistable pairs with
fixed rank $r$ and Hilbert polynomial $rp_E(n)$ is bounded.
\end{proposition}
\begin{proof} We can consider any fixed integer $n\ge 1$ so as all
admissible schemes $\widetilde S$ are immersed into projective
space $\mathbb{P}=P(H^0 (\widetilde S, \widetilde L^n))$ as closed
subschemes and all projective spaces $\mathbb{P}=P(H^0 (\widetilde
S, \widetilde L^n))$ are isomorphic by the construction. Define
$\chi(n)$ as Hilbert polynomial of any subscheme $j: \widetilde S
\hookrightarrow \mathbb{P}$ i.e. $\chi(n):= \chi(j^{\ast}
\mathcal{O}_{\mathbb{P}}(n))$. Then consider corresponding Hilbert
scheme $\mathbb{H}:=\mathrm{Hilb}^{\chi(n)} \mathbb{P}$ and
universal subscheme $\mathbb{U}:=\mathrm{Univ}^{\chi(n)}
\mathbb{P}$. Hilbert scheme $\mathbb{H}$ is (projective) scheme of
finite type over $k$ and the structure morphism $\pi:
 \mathbb{U} \to \mathbb{H}$ is (flat)
morphism of finite type. Hence the scheme $ \mathbb{U}$ is of
finite type over $k$.

Then we can consider families of semistable (not necessarily
locally free) coherent sheaves  on fibres of the morphism $\pi$.
Semistability is understood in usual sense due to Gieseker,
without requirement of quasi-ideality. Sheaves are considered to
have Hilbert polynomial equal to $rp_E(n)$ when it is compute with
respect to the polarization $\widetilde L$. By theorem \ref{bndth}
the set of semistable coherent sheaves with fixed Hilbert
polynomial on fibres of the morphism $\pi$ is bounded. This set is
equal to the set of all pairs $((\widetilde S, \widetilde L),
\widetilde E)$ where $(\widetilde S, \widetilde L)$ be a
2-dimensional projective polarized scheme such that its
polarization $\widetilde L$ induces closed immersion $j:
\widetilde S \hookrightarrow \P$ of $\widetilde S$ as a subscheme
with Hilbert polynomial $\chi(n)$. $\widetilde E$ is
Gieseker-semistable coherent sheaf with Hilbert polynomial equal
to $rp_E(n)$.

 Since admissible semistable
pairs form a subset in the set of all  pairs with
Gieseker-semistable coherent sheaves $\widetilde E$, this
completes the proof of boundedness for families of admissible
semistable pairs.\end{proof}

By  proposition \ref{bnd} there is (common for all semistable
torsion-free sheaves $\widetilde E$) integer $m\gg 0$ such that
the sheaves $\widetilde E \otimes \widetilde L^m$ are globally
generated and for all $\widetilde E$ the vector spaces of global
sections $H^0(\widetilde S, \widetilde E\otimes \widetilde L^m)$
are isomorphic to $k$-vector space $V$ of dimension $rp_E(m)$.
This $V$ will be fixed from now.

 Now we can turn to the  Grassmannian
variety $G(V,r)$ to perform what was done in \cite{Tim4} for
reduced case. Each pair $((\widetilde S,\widetilde L), \widetilde
E)$ with semistable locally free sheaf $\widetilde E$ provides a
closed  immersion $j: \widetilde S \hookrightarrow G(V,r)$. Let
$\OO_{G(V,r)}(1)$ is positive generator of the group $\Pic
G(V,r)$, $P(n):= \chi(j^{\ast}\OO_{G(V,r)}(n))$ be the Hilbert
polynomial of closed subscheme $j(\widetilde S)$. We form Hilbert
scheme $\Hilb^{P(n)}G(V,r)$ and subscheme $H$ formed by all
admissible semistable pairs. This subscheme is defined by
scheme-theoretic images of all possible bases $T$ of families of
S-semistable admissible pairs under induced morphisms to Hilbert
scheme $\Hilb^{P(n)}G(V,r)$. Since we put no restriction to $T$ to
be reduced then it is clear that $H$ can be nonreduced.

%For this purpose regret to the classical construction of Gieseker -- Maruyama moduli scheme by means of GIT-quotient of Quot-scheme. Let $Q$ be the quasiprojective subscheme in $Quot$ and consider its Fitting stratification. The stratum $Q_i\subset Q$ is defined set-theoretically as $Q_i:=\{ q: V \otimes L^{-n} \twoheadrightarrow E| \mathrm{colength} \mathcal{F}itt^0 \mathcal{E}xt^1(E, \mathcal{O}_S) = i\}$ These strata are locally closed in $Q$ and hence quasiprojective.

Now consider a connected component $H_0$ of $H$ containing
$\mu(\widetilde Q)$ mentioned in \cite{Tim4}. Indeed
$\mu(\widetilde Q)$ is the reduction of $H_0$.

\begin{proposition}\label{qpr} {$H_0$ is quasi-projective subscheme in
$\Hilb^{P(t)} G(V,r)$.}
\end{proposition}
\begin{proof}
To prove the quasi-projectivity of $H_0$ we need some topological
result from EGA \cite[ch. 1, proposition 4.5.14]{EGA}. In
particular this
proposition says that the scheme morphism $f: X\to Y$ is %injective
open (resp. closed, homeomorphic onto its image) if and only if
the same holds for its reduction $f_{red}$.

Consider locally closed immersion $\mu(\widetilde Q)
\hookrightarrow \Hilb^{P(n)} G(V, r)$. The Hilbert scheme
$\Hilb^{P(n)} G(V, r)$ is projective $k$-scheme of finite type and
$\mu(\widetilde Q)$ is quasi-projective subscheme in it. The
scheme $H_0$ is another subscheme in $\Hilb^{P(n)} G(V, r)$ such
that $(H_0)_{red}= \mu(\widetilde Q)$. To convince that $H_0$ is
quasi-projective, form scheme-theoretic closures
$\overline{\mu(\widetilde Q)}$ and $\overline H_0$ for
$\mu(\widetilde Q)$ and $H_0$ respectively. Both closures are
taken in $\Hilb^{P(n)} G(V, r)$. Since subschemes
$\overline{\mu(\widetilde Q)}$ and $\overline H_0$ are closed in
the projective scheme they are projective. By the construction
$(\overline H_0)_{red}=\overline{\mu(\widetilde Q)}.$  Now
consider  immersions $i:H_0 \hookrightarrow \overline {H_0}$ and
$i_{red}: \mu(\widetilde Q) \hookrightarrow
\overline{\mu(\widetilde Q)}$. By \cite[ch. 1, proposition
4.5.14]{EGA} since $i_{red}$ is open then $i$ is also open. This
implies that $H_0$ is quasi-projective scheme.\end{proof}

 The group
$PGL(V)$ acts upon the Grassmann variety $G(V, r)$ by linear
transform\-ations of the vector space $V$ and it acts in an
induced fashion upon the product $\Hilb^{P(n)}G(V,r) \times
G(V,r)$. The subscheme $H_0\subset \Hilb^{P(n)}G(V,r)$ is
$PGL(V)$-in\-variant.

\begin{proposition} There exists GIT-quotient $\widetilde M=H_0//PGL(V)$
which is moduli scheme for the functor $\mathfrak f$. The scheme
$\widetilde M$ is a projective algebraic Noetherian scheme of
finite type.
\end{proposition}

\begin{proof}
The scheme $H_0$ is acted upon by the same algebraic group
$PGL(V)$ as its reduction $\mu(\widetilde Q)$ and geometric
invariant theory is also applicable.

Let $\SS$ be the universal quotient bundle on the Grassmannian
$G(V.r)$, as usually $\OO_{G(V,r)}(1)$ is the positive generator
in its Picard group. We use following notations for projections of
the universal subscheme $\Hilb^{P(n)}G(V,r)
\stackrel{\pi}{\longleftarrow } \Univ^{P(n)}G(V,r)
\stackrel{\pi'}{\longrightarrow} G(V,r)$. Form following
 sheaves on the Hilbert scheme $\widetilde L_l^h=\det
\pi_{\ast} \pi'^{\ast} \SS (l)$. Since the projection
 $\pi: \Univ^{P(n)}G(V,r) \to \Hilb^{P(n)}G(V,r)$ is a flat
 morphism and sheaves $\SS(l)$ are locally free, then sheaves
$\widetilde L_l^h$ are invertible.

\begin{proposition} \cite[proposition 19]{Tim4} Sheaves $\widetilde L_l^h$ are very ample for $l\gg 0$.
\end{proposition}
Fix the notation $\widetilde L_l:=\widetilde L_l^h |_{H_0}.$

We remind the following

\begin{definition} \cite[definition 4.2.5]{HL} Let $Y$ be a
$k$-scheme of finite type, $G$ an algebraic $k$-group, and
$\alpha: Y\times G \to Y$ -- group action. {\it $G$-linearization}
of a quasicoherent  $\OO_Y$-sheaf $F$ is an isomorphism of
$\OO_{Y\times G}$-sheaves $\Lambda: \alpha^{\ast}F \to p_1^{\ast}
F$ where $p_1: Y\times G \to Y$ is the projection and the
following cocycle condition holds:
\begin{equation}\label{cocyclel} (\id_Y\times \mult)^{\ast} \Lambda
=p_{12}^{\ast} \Lambda \circ (\alpha \times \id_G)^{\ast} \Lambda.
\end{equation}
Here $p_{12}: Y\times G \times G \to Y \times G$ is a projection
onto first two factors, $\mult: G\times G\to G $ is a morphism of
group multiplication in $G$.
\end{definition}

\begin{proposition} \cite[proposition 20]{Tim4}
Sheaves $\widetilde L_l$ carry $GL(V)$-linear\-ization.
\end{proposition}

Now consider \cite[ch. 4, sect. 4.2]{HL} an arbitrary
one-parameter subgroup $\lambda: \A^1 \setminus 0 \to GL(V).$ We
denote the image of the point $t\in \A^1 \setminus 0$ under the
morphism $\lambda$ by the symbol $\lambda (t)$. The composite of
the morphism
 $\lambda$ with the action $\alpha$ leads to the morphism
  $\alpha(\lambda): \A^1 \setminus 0 \to \mu(\widetilde
Q)$
 for any closed point $\widetilde x\in \mu(\widetilde Q)$.
 This morphism is given by the correspondence  $t\mapsto \widetilde
x_t=\alpha(\lambda(t),\widetilde x)$.

We claim that the morphism $\alpha(\lambda): \A^1 \setminus 0 \to
H_0$ can be uniquely extended to the morphism
$\overline{\alpha(\lambda)}: \A^1 \to H_0$.

Indeed, existence and uniqueness of the extension to the morphism
 $\overline{\alpha(\lambda)}: \A^1 \to \Hilb^{P(n)}G(V,r)$
are provided by properness of the Hilbert scheme. It rests to
convince that the point $\overline{\alpha(\lambda)}(0)$ which is
added to form the closure of the image, belongs to the subscheme
$H_0$.

Turn to  $G=PGL(V)$-equivariant surjective morphisms
\begin{equation*}\xymatrix{\mu(\widetilde Q)&\ar[l]_>>>>>>{\mu'}
 \widetilde Q' \ar[r]^{\phi'}& Q}\end{equation*}
 where $\phi'$ is equivariant standard resolution and $\mu'$ the
 morphism to the Hilbert scheme $\Hilb^{P(n)}G(V,r).$
 
Let $\widetilde x'$ be any point in the preimage
 $\mu'^{-1}(\widetilde x)\subset \widetilde
 Q'$; $x=\phi'(\widetilde x')$. According to
 \cite[lemma 4.3.8]{HL}, the morphism $\beta(\lambda):
 \A^1 \setminus 0 \to Q$ defined by the point $x$ is extended to
 the morphism
 $\overline{\beta(\lambda)}: \A^1 \to Q$. Consider the action
 $\widehat \alpha: G \times \widetilde Q' \to \widetilde
 Q'$ and an induced morphism  $\widehat \alpha_{\lambda}:
 (\A^1 \setminus 0) \times \mu'^{-1}(\widetilde x) \to \widetilde
 Q'$. It includes into the diagram
 \begin{equation*}\xymatrix{(\A^1 \setminus 0) \times
 \mu'^{-1}(\widetilde x) \ar[d]
 \ar[r]^>>>>>>{\widehat \alpha_{\lambda}}& \widetilde Q'
 \ar[d]^{\phi'}\\
 (\A^1 \setminus 0) \times \phi'\mu'^{-1}(\widetilde x)
 \ar[r]^>>>>>>{\beta_{\lambda}}&  Q }
\end{equation*}
where the lower morphism is induced by the action $\beta: G\times
Q \to Q$. The morphism  $\beta_{\lambda}$ extends to the morphism
 $\overline \beta_{\lambda}: \A^1 \times
\phi'\mu'^{-1}(\widetilde x) \to Q$.

Form a Cartesian square
\begin{equation*}
\xymatrix{\A^1  \times \phi'^{-1} (\phi(\mu'^{-1}(\widetilde x)))
\ar[d]_{(\id_G, \phi')}
 \ar[r]^>>>>>>{\widetilde \alpha_{\lambda}}& \widetilde Q'
 \ar[d]^{\phi'}\\
 \A^1  \times \phi'(\mu'^{-1}(\widetilde x))
 \ar[r]^>>>>>>{\overline \beta_{\lambda}}&  Q }\end{equation*}
 where $\widetilde \alpha_{\lambda}$ is a new morphism
 which is defined by equivariance of the morphism $\phi'$ and includes
 in the commutative diagram
 \begin{equation*} \xymatrix{(\A^1 \setminus 0) \times
 \mu'^{-1}(\widetilde x)
 \ar@{^(->}[d] \ar[r]^>>>>>>{\widehat \alpha_{\lambda}}&
 \widetilde Q' \ar@{=}[d]\\
 \A^1 \times \phi'^{-1} (\phi'(\mu'^{-1}(\widetilde x)))
 \ar[r]^>>>>>>{\widetilde \alpha_{\lambda}}& \widetilde Q'}
 \end{equation*}
The immersion in this diagram is induced by immersions of factors.

Existence of the morphism $\widetilde \alpha_{\lambda}$ points, in
particular, on the following fact. The morphism $\widehat
\alpha(\lambda): \A^1 \setminus 0\to \widetilde Q'$ defined by any
point  $\widetilde x' \in \mu'^{-1}(\widetilde x)$ can be extended
to the morphism $\overline{\widehat \alpha(\lambda)}: \A^1 \to
\widetilde Q'$. Descending relatively to the surjective morphism
 $\mu'$ we come to  existence of the extension $\overline{\alpha(\lambda)}:
\A^1 \to \mu(\widetilde Q)$ for the morphism $\alpha(\lambda):
\A^1 \setminus 0 \to \mu(\widetilde Q)$ defined by an arbitrary
point  $\widetilde x \in \mu(\widetilde Q).$

Since  $(H_0)_{red}=\mu(\widetilde Q)$, the morphism
$\alpha(\lambda): \A^1 \setminus 0 \to H_0$ has a unique extension
$\overline{\alpha(\lambda)}: \A^1 \to H_0$.

 Then the point $\widetilde
x_0=\overline{\alpha(\lambda)}(0)$ is a fixpoint of the action of
the subgroup $\lambda$. Notation: $\widetilde x_0=\lim_{t\to 0}
\lambda(t)(\widetilde x).$ The subgroup  $\lambda$ acts on the
fibre $L_{\widetilde x_0}$ of $G$-linearized vector bundle
 $L$ with some weight $r$. Namely, if $\Lambda$ is the
 linearization on $L$ then  $\Lambda (\widetilde x_0, g)=
g^r \cdot \id_{L_{\widetilde x_0}}.$ Define the {\it weight} of
the corresponding one-dimensional representation of the group
$\lambda$ as $w^{\widetilde L_l}(\widetilde x, \lambda)=-r.$

The main tool to analyze the existence of a group quotient is
numerical Hilbert -- Mumford criterion. Recall
\begin{definition} \cite[Definition 4.2.9]{HL} The point  $x\in X$
of the scheme $X$ is {\it semistable with respect to
$G$-linearized ample vector bundle $L$} if there exist an integer
$n$ and an invariant global section $s\in H^0(X, L^{n})$ such that
$s(x)\ne 0$. The point $x$ is {\it stable} if in addition the
stabilizer $\Stab(x)$ is finite and $G$-orbit of the point $x$ is
closed in the open set of all semistable points in $X$.
\end{definition}
\begin{theorem} {\bf (Hilbert -- Mumford criterion)} \rm{ \cite[Theorem
4.2.11]{HL}} {\it The point  $x\in X$ is semistable if and only if
for all nontrivial one-parameter subgroups $\lambda: \A^1\setminus
0\to G$ there is a following inequality $$ w(x,\lambda)\ge 0.
$$
The point $x$ is stable if and only if for all $\lambda$ strict
inequality holds.}
\end{theorem}

Hilbert -- Mumford criterion operates with the set of closed
points of the scheme of interest. The set of closed points of
$H_0$ is the same as the set of closed points in $\mu(\widetilde
Q)$. Hence we can transfer the result of \S 9 \cite{Tim4} into our
situation.

\begin{theorem}\label{quotheor}\cite[theorem 4.2.10]{HL} {Let $G$ be a reductive group
acting on a project\-ive scheme $Y$ with a $G$-linearized ample
line bundle $L$. Then there is a projective scheme $X$ and a
morphism $f: Y^{ss}(L)\to X$ such that $f$ is a universal good
quotient for the $G$-action. Moreover, there is an open subset
$X^s \subset X$ such that $Y^s(L)=f^{-1}(X^s)$ and such that $f:
Y^s \to X^s$ is a universal geometric quotient. }\end{theorem}

We apply this theorem in the following situation:
$X=\overline{H_0},$ $G=PGL(V),$ $L= \widetilde L_l$, $l\gg 0$.
Since we do not know if the equality $(\overline {H_0})^{ss}=H_0$
holds,  we have the following proposition.
\begin{proposition} There is a quasiprojective algebraic scheme
$\widetilde M$ with a morphism  $f: H_0^{ss} \to \widetilde M$,
and $f$ is a universal good $PGL(V)$-quotient. The scheme
$\widetilde M$ contains an open subscheme $\widetilde M^s \supset
\widetilde M$ such that the restriction $f|_{H_0^s}: H_0^s \to
\widetilde M$ is a universal geometric quotient.
\end{proposition}
\begin{proof} Indeed the boundary $\partial \overline H_0= \overline H_0 \setminus
H_0$ is $PGL(V)$-invariant subscheme. It is closed in the
projective scheme $\overline H_0$. This implies that $\partial
\overline H_0$ is projective itself. By Theorem \ref{quotheor} the
formation of $PGL(V)$-quotients leads to the immersion of
projective subscheme $(\partial \overline H_0)^{ss}/PGL(V)
\hookrightarrow (\overline H_0)^{ss}/PGL(V)$ into projective
scheme. Then $\widetilde M=((\overline H_0)^{ss}/PGL(V))\setminus
((\partial \overline H_0)^{ss}/PGL(V))$ what implies
quasi-projectivity of $\widetilde M$. The rest of the proposition
follows immediately from Theorem \ref{quotheor}.
\end{proof}

The same is true for the proof that GIT-quotient is indeed moduli
space as required: the reasoning of \S 12 \cite{Tim4}  is
transferred in our situation directly.

The GIT-quotient $\widetilde M= H_0^{ss} // PGL(V)$ is
quasi-projective Noetherian scheme of finite type over $k$. This
means that there is a locally closed immersion of the scheme
$\widetilde M$ into some appropriate projective $k$-space $\P$.
Let $j: \widetilde M \hookrightarrow \P$ be the corresponding
scheme morphism. By \cite[chapitre 1, proposition 4.5.15]{EGA} if
$f: X\to Y$ is open (resp. closed) immersion the same is true for
$f_{red}$. Then the reduction $j_{red}: \widetilde M_{red} \to \P$
is also locally closed immersion. Now we need an obvious lemma.
\begin{lemma}\label{l} Let $X$ be a projective scheme and
$j: X\hookrightarrow \P$ its locally closed immersion into some
projective space. Then $j$ is closed immersion.
\end{lemma}

The application of lemma \ref{l} yields that $\widetilde M$ is a
projective scheme. The finiteness of the type and Noetherian
property descend from the corresponding properties of the scheme
$H_0$ \cite[ch. 1, \S 2, theorem 1.1]{MumFo}. \end{proof}
\begin{proof}[Proof of the lemma \ref{l}] Since $X$ is projective then there is a projective
space $\P'$ together with closed immersion $j': X\hookrightarrow
\P'$. Two immersion morphisms $j,j'$ induce the morphism $(j',j):
X\hookrightarrow \P' \times \P$ of locally closed immersion. We
have a commutative diagram
\begin{equation*}\xymatrix{X \ar@{_(->}[d]_{j'} \ar@{_(->}[rrd]^{(j',j)}
\ar@{^(->}[rr]^j&& \P\\
\P'&&\ar[ll]^{p_1}\P' \times \P \ar[u]_{p_2}}
\end{equation*}
where $p_1, p_2$ are projections onto the first and the second
factor respectively. Repre\-sentation $j'=p_1 \circ (j', j)$
guarantees that $(j',j)$ is closed immersion \cite[ch. 1, remarque
5.2.7 (iii)]{EGA}. Then by the representation $j=p_2 \circ (j',j)$
we conclude that $j$ is closed immersion.
\end{proof}

\begin{remark} Taking into account the result of \S 10 of the article
\cite{Tim4} we come to the following conclusion. The morphism of
compactifications constructed there holds when $\widetilde M$ and
$\overline M$ are reduced or when nonreduced moduli schemes
$\overline M$ and $\widetilde M$ are replaced by their reductions.
The existence of such a morphism for nonreduced schemes is open
question for now. \end{remark}

\section{Subfunctors and moduli subschemes}

Let $\mathfrak f: (Schemes_k)^o \to (Sets)$  be a functor
assigning to a scheme  $T$ the set of families $\mathfrak F_T$
with base $T$. Families are mentioned to consist of objects of
some prescribed class $\mathfrak F$. Let $P$ be some property,
$\mathfrak F^P$ subclass of objects in $\mathfrak F$ with the
property $P$. $P$ is called {\it open} if for any family
$\mathfrak F_T$ points  $t\in T$ corresponding to objects with
$P$, constitute an open subscheme $T_0$ in $T$. Then it makes
sense to call the functor  $\mathfrak f^P: (Schemes_k)^o \to
(Sets)$ assigning to a scheme  $T$ the set of families $\mathfrak
F^P_T$ of objects from the class $\mathfrak F^P$, as {\it open
$P$-subfunctor} of the functor  $\mathfrak f$.

\begin{proposition}\label{submoduli} If $\mathfrak f$ is corepresented by
the scheme $M$ then its open subfunctor $\mathfrak f^P$ is
corepresented by an open subscheme  $M_0$ in $M$.
\end{proposition}
\begin{proof} Corepresentability of the functor  $\mathfrak f$ by the scheme
$M$ means in particular the commutativity of diagrams of the view
\begin{equation*} \xymatrix{T \ar@{|->}[d]_{\underline M}
\ar@{|->}[r]^{\mathfrak f}&
\ar[ld]_{\psi_T} \{\mathfrak F_T\} \ar[d]^{\psi'_T} \\
\Hom(T,M) \ar[r]^{\omega_T}& \Hom (T,F')}
\end{equation*}

Let  $M_0$ be the union of images
$$\bigcup_{\small {\begin{array}{c} f\in \psi_T(\{\mathfrak
F^P_T\})\\ T\in Ob (Schemes_k)
\end{array}}} f(T)$$ in $M$.

In the sequel when working with the subclass $\mathfrak F^P$ and
with the subfunctor $\mathfrak f^P$ we use superscript $P$ in
notations of natural transformations and induced morphisms in the
category of sets: $\psi^P,$ $\psi^P_T$  etc.

It follows from the construction that the subfunctor $\mathfrak
f^P$ is corepresented by the subscheme $M_0$. To prove this it is
enough to construct a natural transformation  $\omega^P$. For the
scheme immersion $M_0 \hookrightarrow M$, for each scheme  $T$ and
for the induced inclusion of sets  $\Hom(T, M_0) \hookrightarrow
\Hom(T, M)$ consider the diagram
\begin{equation*}\xymatrix{\{\mathfrak F^P_T\} \ar[ddr]_{\psi_T}
\ar[dr]^>>>{\psi^P_T} \ar[drr]^>>>>>>>>{\psi'_T={\psi'}^P_T}\\
& \Hom(T, M_0) \ar@{^(->}[d] \ar@{.>}[r]_{\omega^P_T}& \Hom(T,
F')\\
&\Hom(T,M) \ar[ur]_{\omega_T } }\end{equation*} It shows that the
map $\omega^P_T$ is uniquely defined as the composite in the lower
triangle. Respectively, for the natural transformation of our
interest we have the definition by the composite
$$\omega^P: \Hom(-, M_0) \to \Hom(-, M)
\stackrel{\omega}{\longrightarrow} \Hom (-, F').$$

It rests to convince that  $M_0$ is  open subscheme in $M$. For
this consider the covering $M=\bigcup _{U\in Ob\,Open(M)}U$ of the
scheme $M$ by all possible open subschemes, and corresponding sets
 $\{\mathfrak F_U\}$. Let  $U_0\subset U$ be the maximal open
 subscheme in $U$ such that the family $\mathfrak F_U$ consists of objects
 with the property $P$ if restricted on $U_0$. Such subscheme
 depends on the choice of family in the set $\{\mathfrak F_U\}$.
 We consider all such subschemes  $U_0$ for each  $U$. It is clear that
 $\bigcup U_0 \subset M_0$ and $\bigcup
U_0$ is an open subscheme in  $M$.

Now let  $T_0$ be an arbitrary scheme and let it be the base of a
family $\mathfrak F^P_{T_0}$. The natural transformation
 $\psi$ assigns to it the morphism  $f: T_0 \to M$. This morphism
 factors through the subscheme  $M_0$. Then form  open covering
 of the scheme $f(T_0)$ by subschemes $U \times_M
T_0$. Images of these subschemes in  $M$ belong to $\bigcup U_0$.
The consideration of different schemes $T_0$ yields in the
equality $M_0=\bigcup U_0$.
\end{proof}

Now turn to the class of objects described in (\ref{class}), and
to the property  $P=\{\pi^{-1}(t)\cong S\}$. It is clear that $P$
is open property and hence there is open subfunctor $\mathfrak
f^P$. The proposition  \ref{submoduli} leads to the existence of
the open subscheme $\widetilde M_0$ which is the coarse moduli
space for the subfunctor $\mathfrak f^P$, i.e. the {\it moduli
scheme of  $S$-pairs}.

Analogously, consider the Gieseker -- Maruyama functor
$\overline{\mathfrak f}: (Schemes_k)^o \to Sets$ assigning to any
scheme  $T$ the set of flat $T$-based families  $\E$ of semistable
coherent torsion-free $\OO_S$-sheaves  $E=\E|_{t\times S}$ with
Hilbert polynomial $\chi (E \otimes L^n)=rp_E(n)$. The property
 $P'=\{E \mbox{ is locally free} \}$ provides the subscheme
 $\overline M_0$ of moduli of semistable vector bundles. This
 subscheme corepresents  the subfunctor
  $\overline {\mathfrak f}^{P'}$.

By definitions of subfunctors $\mathfrak f^P$ and
$\overline{\mathfrak f}^{P'}$ they are isomorphic: corresponding
classes of objects of parametrization coincide. This implies the
isomorphism of moduli schemes  $\widetilde M_0 \cong \overline
M_0$. This completes the proof of the theorem \ref{th1}.

{\bf Acknowledgement.} The author expresses her deep and sincere
gratitude to Prof. Dr. Vasile Brinzanescu (IMAR, Bucharest,
Romania) for interest to this work and fruitful discussions. Also
author thanks the Institute of Mathematics of Romanian Academy
(IMAR) where the most part of the work was done, for hospitality
and support.

\end{document}